\documentclass[reqno]{amsart}
\usepackage[english]{babel}
\usepackage{amssymb}
\usepackage[latin1]{inputenc}
\newtheorem{theorem}{Theorem}[section]

\newtheorem{lemma}[theorem]{Lemma}

\theoremstyle{definition}

\theoremstyle{remark}

\numberwithin{equation}{section}

\newcommand{\CEP}[1]{\mbox{$\mathbb{C}^{#1}$}}

\newcommand{\supp}{\mbox{supp}}
\newcommand{\E}{\mbox{$\mathcal{E}$}}

\newcommand{\Eo}{\mbox{$\mathcal{E}_{0}$}}

\newcommand{\Ep}{\mbox{$\mathcal{E}_{p}$}}
\newcommand{\Ec}{\mbox{$\mathcal{E}_{\chi}$}}

\newcommand{\F}{\mbox{$\mathcal{F}$}}

\newcommand{\ddcn}[1]{\mbox{$\left(dd^{c}#1\right)^{n}$}}

\newcommand{\n}{\mbox{$\mathcal{N}$}}

\begin{document}
\title{On a  Monge-Amp\`{e}re type equation in the Cegrell class $\mathcal{E}_{\chi}$}

\author{Rafa\l\ Czy{\.z}}
\address{Jagiellonian University\\ \L ojasiewicza 6\\ 30-348 Krak\'ow\\ Poland}
\email{Rafal.Czyz@im.uj.edu.pl}
\thanks{The author was partially supported by ministerial grant number N N201 3679 33}

\keywords{complex Monge-Amp\`{e}re operator, Cegrell classes, Dirichlet problem, plurisubharmonic function}
\subjclass[2000]{Primary 32W20; Secondary 32U15.}
%\date{\today}
%
\begin{abstract}
We prove an existence and uniqueness  result for a Monge-Amp\`{e}re
type equation in the Cegrell class $\E_{\chi}$.
\end{abstract}

\maketitle

\section{Introduction}

 It is a classical problem in analysis to find, for given function $F$,
solutions $u$ to the equation:
\begin{equation}\label{intr}
(dd^cu)^n=F(z,u(z))d\mu\, ,
\end{equation}
where $(dd^cu)^n$ is the complex Monge-Amp\`{e}re operator. Equations of the type~(\ref{intr}) has
played a significant not only within the fields of fully nonlinear second order elliptic equations
and pluripotential theory, but also in applications. We refer
to~\cite{bt,bt2,ckns,c,cegrell_kolodziej,kolweak,krylov} and the reference therein for further
information about equations of Monge-Amp\`{e}re type.

 In the last decade the so called Cegrell classes has played a prominent role in working
with the complex Monge-Amp\`{e}re operator. For further
information about the Cegrell classes see
e.g.~\cite{ahag_thesis,ahag_Ff,czyz_compl,czyz_cegrell,
czyz_energy,cegrell_energy,cegrell_gdm,cegrell_bdd}.
Let $\Eo$, $\Ep$, $\F$, $\n$,
and $\E$ be as in~\cite{cegrell_energy,cegrell_gdm,cegrell_bdd}.

In~\cite{bgz}, Benelkourchi, Guedj, and Zeriahi introduced the
following notation of the Cegrell classes. Let
$\chi:(-\infty,0]\to (-\infty,0]$ be a continuous and
nondecreasing function and let $\Ec$ contain those
plurisubharmonic functions $u$ for which there exists a decreasing
sequence $u_{j}\in\Eo$, which converges pointwise to $u$ on
$\Omega$, as $j$ tends to $+\infty$, and
\[
\sup_{j}\int_{\Omega} -\chi(u_j)\ddcn{u_j}< \infty\, .
\]
Note that if $\chi=-(-t)^p$, then $\Ec=\Ep$, and if $\chi=-1$,
then $\Ec=\F$.

 The measure $\ddcn{u}$ might have infinite total mass, i.e.
$(dd^cu)^n(\Omega)=+\infty$. On the other hand, if $u\in
\Ec(\Omega)$, then under certain assumption on the function
$\chi$, the measure $(dd^cu)^n$ vanishes on all pluripolar sets in
$\Omega$ and the following integral is always finite
\[
\int_{\Omega}-\chi(u)(dd^cu)^n<+\infty,
\]
which means that $-\chi(u)(dd^cu)^n$ is a positive and finite
measure on $\Omega$. For this reason it is natural to consider the
following Monge-Amp\`{e}re type equation:
\[
-\chi(u)(dd^cu)^n=d\mu\, ,
\]
where $\mu$ is a given positive and finite measure on $\Omega$,
which vanishes on pluripolar sets in $\Omega$. We prove the
following theorem.

\bigskip

\noindent {\bf Main Theorem.} {\em Assume that $\Omega$ is a bounded hyperconvex domain in
$\CEP{n}$, $n\geq 1$, and let $\chi:(-\infty,0]\to (-\infty,0]$ be a continuous increasing function
such that $\chi(0)=0$ and $\lim_{t\to -\infty}\chi(t)=-\infty$, such that $\Ec\subset \E$. If $\mu$
is a positive and finite measure in $\Omega$, such that $\mu(P)=0$, for all pluripolar sets
$P\subset \Omega$, then there exists a function $u\in \Ec$ such that}
\[
-\chi(u)(dd^cu)^n=d\mu\, .
\]
{\it Furthermore, if $\Ec\subset \n$, then the solution of above equation is
uniquely determined}.

\bigskip
$\E$ denotes the Cegrell class of all negative plurisubharmonic functions on which the complex
Monge-Amp\`{e}re operator $(dd^c\cdot)^n$ is well defined and $\n\subset \E$ denotes the Cegrell
class for which the smallest maximal plurisubharmonic majorant is identically equal to $0$. It
follows from~\cite{cegrell_energy}, \cite{cegrell_gdm}, \cite{cegrell_bdd} that $\Ep,
\F\subseteq\n$.

It was proved in~\cite{b} that if $\chi:(-\infty,0]\to (-\infty,0]$ is a continuous increasing,
convex or concave function such that $\chi(0)=0$ and $\lim_{t\to -\infty}\chi(t)=-\infty$, then
$\Ec\subset \E$.

It should be noted that uniqueness part is an immediate
consequence of Cegrell's work in~\cite{cegrell_bdd} and the fact
that $\Ec$ is properly contained in $\n$. Let $\chi(t)=-(-t)^p$,
$(p>0)$, in the Main Theorem, then it follows that: If $\mu$ is a
positive and finite measure in $\Omega$, such that $\mu(P)=0$, for
all pluripolar sets $P\subset \Omega$, then there exists a unique
function $u\in \Ep$ such that
\[
(-u)^p(dd^cu)^n=d\mu\, .
\]
If $\chi:(-\infty,0]\to (-\infty,0]$ is a continuous function such
that $\chi(0)<0$ and $\lim_{t\to -\infty}\chi(t)=-\infty$, then
the existence of solution to the Monge-Amp\`{e}re type equation
given by
\[
-\chi(u)(dd^cu)^n=d\mu\,
\]
is a consequence of~\cite{cegrell_kolodziej} with the assumption that $-\chi(t)^{-1}$ is bounded.

\bigskip

The author would like to thank Per \AA hag and S{\l}awomir
Ko{\l}odziej for their generous help and encouragement. This
research was done in part during the authors' visit to the
University of Vienna in 2007. He would like to thank the members
of the department of Mathematics for their kind hospitality.

\section{Proof of the Main Theorem}

\begin{lemma}\label{lep}
Let $\Omega$ be a bounded hyperconvex domain in $\CEP{n}$. If a
sequence $u_j\in\F$ satisfies the condition
\[
\sup_j\int_{\Omega}(dd^cu_j)^n<+\infty,
\]
and if there exists $u\in PSH(\Omega)$ such that $u_j\to u$
weakly, then $u\in \F$.
\end{lemma}

\begin{proof}
From~\cite{cegrell_gdm} there exists $w_j\in \Eo\cap \mathcal C(\bar \Omega)$ such that $w_j\searrow
u$, $j\to \infty$. Note that since $u_j\to u$ weakly, then $u=\lim_{j\to \infty}v_j$, where
\[
v_j=\left(\sup_{k\geq j}u_k\right)^*.
\]
Observe that $v_j$ is a decreasing sequence, $v_j\geq u_j$, so $v_j\in \F$ and
from~\cite{cegrell_gdm} we have
\[
\int_{\Omega}(dd^cv_j)^n\leq \int_{\Omega}(dd^cu_j)^n.
\]
Define
\[
\varphi_j=\max(w_j,v_j).
\]
Then $\varphi_j\in \Eo$, $\varphi_j$ is a decreasing sequence $v_j\searrow u$ and again
from~\cite{cegrell_gdm} we get
\[
\sup_j\int_{\Omega}(dd^c\varphi_j)^n\leq \sup_j\int_{\Omega}(dd^cv_j)^n\leq
\sup_j\int_{\Omega}(dd^cu_j)^n<+\infty,
\]
which means that $u\in \F$.
\end{proof}

First we prove our Main Theorem in the case of compactly supported measures.

\begin{lemma}\label{lemma}
Let $\Omega$ be a bounded hyperconvex domain in $\CEP{n}$, and let
$\chi:(-\infty,0]\to (-\infty,0]$ be a continuous increasing
function such that $\chi(0)=0$ and $\lim_{t\to
-\infty}\chi(t)=-\infty$. If $\mu$ is a positive, finite, and
compactly supported measure in $\Omega$, such that $\mu(P)=0$ for
all pluripolar sets $P\subset \Omega$, then there exists a unique
function $u\in \F\cap\Ec$ such that
\begin{equation}\label{mat}
-\chi(u)(dd^cu)^n=d\mu.
\end{equation}
\end{lemma}

\begin{proof}
If $\mu\equiv 0$, then it is clear that $u=0$ is a solution of
(\ref{mat}). Assume now that $\mu\not\equiv 0$. For $k\in \Bbb N$
consider the following equations
\begin{equation}\label{ck}
(dd^cu_k)^n=\min\left(\frac {-1}{\chi (u_k)},k\right)d\mu\, .
\end{equation}
The function defined by
\[
F_k(t)=\min\left(\frac {-1}{\chi (t)},k\right)
\]
is bounded and continuous. Therefore it follows from~\cite{cegrell_kolodziej} that
there exists $u_k\in \F$, which satisfies (\ref{ck}). We also have that
\[
(dd^cu_k)^n=\min\left(\frac {-1}{\chi (u_k)},k\right)d\mu\leq \frac {-1}{\chi
(u_k)}d\mu,
\]
so
\[
\sup_k\int_{\Omega}(-\chi (u_k))(dd^cu_k)^n\leq
\mu(\Omega)<+\infty.
\]
We shall next prove that there exist $\alpha\in \Eo,\beta\in \F$
such that
\begin{equation}\label{bineq}
\beta \leq u_k\leq \alpha
\end{equation}
almost everywhere $[d\mu]$, for $k\geq 2$.

By Cegrell decomposition theorem (see~\cite{cegrell_gdm}) there
exist $\phi \in \Eo$ and $f\in L^1((dd^c\phi)^n)$, $f\geq 0$ such
that
\[
\mu=f(dd^c\phi)^n.
\]
Fix $a>0$ such that $\chi (-a)\geq -\frac 12$. Then there exists
$\alpha\in \Eo$ such that (see~\cite{kol})
\[
(dd^c\alpha)^n=\min\left(f,\frac {a^n}{\|\phi\|^n}\right)(dd^c\phi)^n,
\]
where $\|\phi\|=\sup_{z\in \Omega}|\phi(z)|$. The comparison principle (see~\cite{bt2}) yields that
\[
\alpha \geq \frac {a}{\|\phi\|}\phi\geq -a
\]
and
\[
\int_{\{\alpha <u_k\} }(dd^cu_k)^n\leq \int_{\{\alpha <u_k\}
}(dd^c\alpha)^n\leq \int_{\{\alpha <u_k\} }d\mu.
\]
Observe that on the set $\{\alpha <u_k\}$ we have $u_k>-a$ and
\[
(dd^cu_k)^n=\min\left(\frac {-1}{\chi (u_k)},k\right)d\mu\geq \min\left(\frac
{-1}{\chi (-a)},k\right)d\mu\geq 2d\mu,
\]
for $k\geq 2$, which implies that $\mu(\{\alpha <u_k\})=0$, $k\geq
2$.

There exists $\psi\in \F$ such that $(dd^c\psi )^n=d\mu$
(see~\cite{cegrell_gdm}). Fix $w \in \Eo$ and $b>0$ such that
\[
\chi (\sup_{supp\, \mu}(\psi+bw))<-2.
\]
Let $\beta=\psi +bw$. Note that $(dd^c\beta)^n\geq d\mu$. By the comparison principle
(see~\cite{bt2}) we obtain
\[
\int_{\{u_k<\beta\} }d\mu\leq \int_{\{u_k<\beta\}
}(dd^c\beta)^n\leq \int_{\{u_k<\beta\} }(dd^cu_k)^n,
\]
but on the set $\{u_k<\beta\}\cap supp\,\mu$ we have
$u_k<\beta\leq \sup_{supp\, \mu}\beta$ and
\[
(dd^cu_k)^n=\min\left(\frac {-1}{\chi (u_k)},k\right)d\mu\leq \frac 12d\mu,
\]
which means that $\mu(\{u_k<\beta\})=0$, for all $k$.

Now it follows from (\ref{bineq}) that there exists a
plurisubharmonic function $u\neq 0$ and a subsequence (also
denoted by $u_k$) such that $u_k\to u$ almost everywhere $[d\mu]$.
Since $u\neq 0$ then
\[
-\frac 1{\chi(\sup_{\text{supp} \,\mu}u)}<+\infty.
\]
By Hartog's lemma, functions
\[
F_k(u_k)=\min(-\chi(u_k)^{-1},k)
\]
are uniformly bounded on $\supp\,\mu$ and therefore
\[
\sup_k\int_{\Omega}(dd^cu_k)^n\leq \sup_k
\int_{\Omega}F_k(u_k)d\mu<+\infty.
\]
Lemma~\ref{lep} yields that $u\in \F$.

The stability theorem proved in~\cite{cegrell_kolodziej} implies
that the weak convergence, $u_k\to u$, is equivalent to
convergence in capacity. Using Xing's theorem in~\cite{xing} we
get that $(dd^cu_k)^n\to (dd^cu)^n$ weakly. Therefore using the
dominated convergence theorem we get that
\[
(dd^cu)^n=\lim_{k\to +\infty}(dd^cu_k)^n=\lim_{k\to
+\infty}F_k(u_k)d\mu=\frac {-1}{\chi(u)}d\mu.
\]
So we have proved that there exists  a solution $u\in \F$ to
(\ref{mat}). Then we have that
\[
\int_{\Omega}(-\chi(u))(dd^cu)^n<+\infty
\]
and therefore it follows from~\cite{bgz} that $u\in \Ec$.

It will be proved in the proof of the Main Theorem that if $u,v\in
\F$ are solutions of (\ref{mat}) then $(dd^cu)^n=(dd^cv)^n$ and
therefore $u=v$ (see~\cite{cegrell_gdm}).
\end{proof}

\begin{proof}[Proof of the Main Theorem] Assume that $\mu$ is a positive and
finite measure in $\Omega$ such that $\mu(P)=0$ for all pluripolar sets $P\subset \Omega$. Let
$\Omega_j$ be a fundamental sequence of strictly pseudoconvex domains, i.e. $\Omega_j\Subset
\Omega_{j+1}\Subset \Omega$ and $\bigcup_{j=1}^{\infty}\Omega_j=\Omega$ (see~\cite{cegrell_bdd}).
Let us define $d\mu_j={\bf 1}_{\Omega_j}d\mu$, where ${\bf 1}_{\Omega_j}$ is a characteristic
function for $\Omega_j$. By Lemma~\ref{lemma} there exists a sequence $u_j\in \F\cap\Ec$ such that
\[
-\chi(u_j)(dd^cu_j)^n=d\mu_j.
\]
We now shall prove that $u_j$ is a decreasing sequence. Let $A=\{z\in
\Omega:u_j(z)<u_{j+1}(z)\}$. On the set $A$, we have that
\[
(dd^cu_j)^n=-\chi(u_j)^{-1}d\mu_j\leq
-\chi(u_{j+1})^{-1}d\mu_j\leq
-\chi(u_{j+1})^{-1}d\mu_{j+1}=(dd^cu_{j+1})^n
\]
and by the comparison principle (see~\cite{bt2}) we get that
\[
\int_A(dd^cu_{j+1})^n\leq \int_A(dd^cu_{j})^n.
\]
Hence,
\begin{equation}\label{1}
(dd^cu_{j})^n=(dd^cu_{j+1})^n
\end{equation}
on $A$. Similarly on the set $B=\Omega_j\setminus A=\{z\in
\Omega_j:u_j(z)\geq u_{j+1}(z)\}$ we obtain that
\begin{equation}\label {2}
(dd^cu_j)^n=-\chi(u_j)^{-1}d\mu_j\geq -\chi(u_{j+1})^{-1}d\mu_j=
-\chi(u_{j+1})^{-1}d\mu_{j+1}=(dd^cu_{j+1})^n.
\end{equation}
From the equalities (\ref{1}) and (\ref{2}) we get that
$(dd^cu_{j})^n\geq(dd^cu_{j+1})^n$ on $\Omega_j$. This implies
that $-\chi(u_j)^{-1}d\mu_j\geq-\chi(u_{j+1})^{-1}d\mu_j$ and then
$\chi(u_j)\geq \chi(u_{j+1})$ a.e. $[d\mu_j]$, so $u_j\geq
u_{j+1}$ a.e. $[d\mu_j]$. Hence $\mu_j(\{u_j<u_{j+1}\})=0$ and
$(dd^cu_j)^n=0$ on $A\cap\Omega_j$. Since $(dd^cu_j)^n=d\mu_j=0$
on $\Omega\setminus \Omega_j$ we finally obtain that
$(dd^cu_j)^n=0$ on $A=\{u_j<u_{j+1}\}$.

Now take $\psi\in \Eo$ such that $(dd^c\psi)^n=d\lambda$, where $d\lambda$ is the Lebesgue measure,
and consider $A_k=\{z\in\Omega:u_j<u_{j+1}+\frac 1k\psi \}$. Observe that $u_{j+1}+\frac
1k\psi\in\F$ and $A_k\subset A$. By the comparison principle (see~\cite{bt2}) we obtain that
\[
\int_{A_k}\ddcn{(u_{j+1}+\frac 1k\psi)}\leq
\int_{A_k}(dd^cu_j)^n\leq \int_{A}(dd^cu_j)^n=0,
\]
and then
\[
0=\int_{A_k}\ddcn{(u_{j+1}+\frac 1k\psi)}\geq \frac
1k\int_{A_k}(dd^c\psi)^n=\frac 1k\lambda(A_k),
\]
which means that $\lambda(A_k)=0$. Hence
$\lambda(A)=0$, since $A=\bigcup_{k=1}^{\infty}A_k$. We have proved
that $u_j\geq u_{j+1}$ a.e $[d\lambda]$, but since functions
$u_j,u_{j+1}$ are plurisubharmonic we obtain that $u_j\geq
u_{j+1}$ on $\Omega$, so $u_j$ is a decreasing sequence.

Note also that
\[
\sup_j\int_{\Omega}-\chi(u_j)(dd^cu_j)^n\leq \int_{\Omega}d\mu<+\infty.
\]
Moreover from~\cite{bgz} we obtain
\[
\int_{\Omega}-\chi(u_j)(dd^cu_j)^n=\int_0^{+\infty}\chi'(t)(dd^cu_j)^n(\{u_j<-t\})dt\geq
\int_0^{+\infty}\chi'(t)t^nC_n(\{u_j<-2t\})dt,
\]
where $C_n$ is Bedford-Taylor capacity (see~\cite{bt2}). Therefore
\[
\sup_j \int_0^{+\infty}\chi'(t)t^nC_n(\{u_j<-t\})dt\leq \sup_j\int_{\Omega}-\chi(u_j)(dd^cu_j)^n\leq
\int_{\Omega}d\mu<+\infty,
\]
which implies that there exists $\lim_{j\to \infty}u_j=u\in\Ec$, and $(dd^cu_j)^n$ tends weakly to
$(dd^cu)^n$. Therefore using the monotone convergence theorem we get that
\[
(dd^cu)^n=\lim_{j\to \infty}(dd^cu_j)^n=\lim_{j\to
\infty}-\chi(u_j)^{-1}{\bf 1}_{\Omega_j} d\mu=-\chi(u)^{-1}d\mu.
\]
This ends the proof of the existence part of this theorem.

Now we will proceed with the uniqueness part. Assume that
$\Ec\subset \n$ and suppose that there exist $u,v\in \Ec$ such
that $-\chi(u)(dd^cu)^n=-\chi(v)(dd^cv)^n=d\mu$. Observe that on
the set $\{z\in \Omega:u(z)<v(z)\}$ we have
\[
(dd^cu)^n=-\chi(u)^{-1}d\mu\leq -\chi(v)^{-1}d\mu=(dd^cv)^n.
\]
Using the comparison principle (see~\cite{bt2}) we obtain
\[
\int_{\{u<v\}}(dd^cv)^n\leq \int_{\{u<v\}}(dd^cu)^n,
\]
so $(dd^cu)^n=(dd^cv)^n$ on $\{z\in \Omega:u(z)<v(z)\}$.
Similarly, we get that $(dd^cu)^n=(dd^cv)^n$ on the set $\{z\in
\Omega:u(z)>v(z)\}$. Since $\mu$ does not put mass on pluripolar
sets and $\{u=-\infty\}=\{\chi(u)=-\infty\}$ and
$\{v=-\infty\}=\{\chi(v)=-\infty\}$, then $(dd^cu)^n=(dd^cv)^n=0$
on the set $C=\{u=-\infty\}\cup\{v=-\infty\}$. On the set
$\{u=v\}\setminus C$ we also have
\[
(dd^cu)^n=-\chi(u)^{-1}d\mu= -\chi(v)^{-1}d\mu=(dd^cv)^n.
\]
Thus, $(dd^cu)^n=(dd^cv)^n$ on $\Omega$, which implies that $u=v$
(see~\cite{cegrell_bdd}).
\end{proof}


\begin{thebibliography}{30}

\bibitem{ahag_thesis} \AA hag P., The complex Monge-Amp\`{e}re operator on bounded
hyperconvex domains, Ph. D. Thesis, Ume\aa \ Univ., 2002.

\bibitem{ahag_Ff} \AA hag P., A Dirichlet Problem for the complex
Monge-Amp\`{e}re Operator in ${\mathcal F}(f)$, Michigan Math. J.
55 (2007), 123-138.

\bibitem{czyz_compl} \AA hag P. \& Czy\.z  R., The connection between the Cegrell
classes and compliant functions, Math. Scand. 99 (2006),
87-98.

\bibitem{czyz_cegrell} \AA hag P. \& Czy\.z  R., On the Cegrell
classes, Math. Z. 256 (2007), 243-264.

\bibitem{czyz_energy} \AA hag P., Czy\.z R. \& Ph\d{a}m H.H., Concerning the energy
class $\Ep$ for $0<p<1$, Ann. Polon. Math. 91 (2007), no. 2-3,
119--130.

\bibitem{bt} Bedford, E., Taylor, B.A., The Dirichlet problem for an
equation of complex Monge-Amp\`{e}re type. In: Byrnes, C. (ed.)
Partial Differential Equations and Geometry. Dekker, 39-50, 1979.

\bibitem{bt2} Bedford E. \& Taylor B.A., A new capacity for
plurisubharmonic functions, Acta Math. 149 (1982), 1-40.

\bibitem{b} Benelkourchi S., Weighted pluricomplex energy, arXiv:0806.4850 (2008).

\bibitem{bgz} Benelkourchi S., Guedj V. \& Zeriahi A.,
Plurisubharmonic functions with weak singularities, Acta Universitatis Upsaliensis, Proceedings of the
conference in honour of C. Kiselman ("Kiselmanfest", Uppsala, May 2006) (in press).

\bibitem{ckns} Caffarelli, L., Kohn, J.J., Nirenberg, L., Spruck, J.,
The Dirichlet problem for nonlinear second order elliptic
equations. II. Complex Monge-Ampère and uniformly elliptic
equations, Comm. of Pure and Appl. Math. 38, (1985), 209-252


\bibitem{c} Cegrell, U., On the Dirichlet problem for the complex Monge-Amp\`{e}re
operator, Math. Z. 185, 247-251, (1984).

\bibitem{cegrell_energy} Cegrell U., Pluricomplex energy, Acta Math.
180, (1998), 187-217.


\bibitem{cegrell_gdm} Cegrell U., The general definition of the complex
Monge-Amp\`{e}re operator, Ann. Inst. Fourier (Grenoble)  54
(2004), 159-179.


\bibitem{cegrell_bdd} Cegrell U., A general Dirichlet
problem for the complex Monge-Amp\`{e}re operator, Ann. Polon.
Math. (in press).


\bibitem{cegrell_kolodziej} Cegrell U. \& Ko\l odziej S., The equations of complex
Monge-Amp\`{e}re type and stability of solutions, Math. Ann.
334 (2006), 713-729.

\bibitem{kol} Ko\l odziej S., The range of the complex Monge-Amp\`{e}re operator,
II, Indiana Univ. Math. J. 44 (1995), 765-782.

\bibitem{kolweak} Ko\l odziej S., Weak solutions of equations of complex Monge-Amp\`{e}re type. Ann.
Polon. Math. 73 (2000), 59-67.

\bibitem{krylov} Krylov N.V., Fully nonlinear second order elliptic
equations: recent development, Ann. Scuola Norm. Sup. Pisa Cl.
Sci. (4) 25 (1997), 569-595 (1998).


\bibitem{xing} Xing Y., Continuity of the complex
Monge-Amp\`{e}re operator, Proc. Amer. Math. Soc., 124 (1996),
457-467.



\end{thebibliography}
\enddocument